\documentclass[11pt,reqno]{amsart}

\usepackage{amsmath,amsthm}
\usepackage{amssymb}
\usepackage{epic,eepic}

\numberwithin{equation}{subsection}

\newcommand{\sqsp}{\renewcommand{\baselinestretch}{1.2}\tiny\normalsize}

\newtheorem{thm}[subsection]{Theorem}
\newtheorem{lemma}[subsection]{Lemma}
\newtheorem{prop}[subsection]{Proposition}
\newtheorem{cor}[subsection]{Corollary}

\theoremstyle{definition}

\newcommand{\biglpren}{\bigl (}  
\newcommand{\bigrpren}{\bigr )}  
\newcommand{\Thetabar}{\overline{\Theta}}
\newcommand{\Thetatilde}{\widetilde{\Theta}}
\newcommand{\thetabar}{\overline{\theta}}

\newcommand{\Lbar}{\overline{L}}
\newcommand{\Rbar}{\overline{R}}
\newcommand{\Mbar}{\overline{M}}

\newcommand{\lprod}{\dashv}  
\newcommand{\rprod}{\vdash}  
\newcommand{\mprod}{\bot}   
\newcommand{\Htrias}{H_{\text{Trias}}}
\newcommand{\HT}{H^{\text{Trias}}}
\newcommand{\Ctrias}{C_{\text{Trias}}}
\newcommand{\CT}{C^{\text{Trias}}}

\DeclareMathOperator{\Id}{Id}
\DeclareMathOperator{\im}{im}
\DeclareMathOperator{\Hom}{Hom}
\DeclareMathOperator{\Ob}{Ob}
\DeclareMathOperator{\Ext}{Ext}
\DeclareMathOperator{\Der}{Der}
\DeclareMathOperator{\Inn}{Inn}

\begin{document}
\title{(Co)homology of triassociative algebras}
\author{Donald Yau}

\begin{abstract}
We study homology and cohomology of triassociative algebras with non-trivial coefficients.  The cohomology theory is applied to study algebraic deformations of triassociative algebras.
\end{abstract}

\email{dyau@math.ohio-state.edu}
\address{Department of Mathematics, The Ohio State University Newark, 1179 University Drive, Newark, OH 43055, USA}

\maketitle
\sqsp


\section{Introduction}
\label{sec:intro}

Triassociative algebras, introduced by Loday and Ronco \cite{lr}, are analogues of associative algebras in which there are three binary operations, satisfying $11$ relations.  These algebras are closely related to some very interesting combinatorial objects.  In fact, the operad $Trias$ for triassociative algebras are modeled by the standard simplices.  This operad is Koszul, in the sense of Ginzburg and Kapranov \cite{gk}.  Its Koszul dual operad $TriDend$ is modeled by the Stasheff polytopes \cite{sta}, or equivalently, planar trees.  The algebras for the latter operad are called tridendriform algebras, which have three binary operations whose sum is an associative operation.

The purpose of this paper is to advance the study of triassociative algebras in the following ways: (1) Algebraic foundation: When studying (co)homology, the first task is to figure out the correct coefficients.  We introduce homology and cohomology with non-trivial coefficients for triassociative algebras.  In fact, Loday and Ronco \cite{lr} already defined homology with trivial coefficients for these algebras.  We will build upon their construction.  (2) Deformation: Once (co)homology is in place, we use it to study algebraic deformations of triassociative algebras, following the pattern established by Gerstenhaber \cite{ger}.

We should point out that deformations of dialgebras \cite{loday}, the analogues of triassociative algebras with only two binary operations, have been studied by Majumdar and Mukherjee \cite{mm}.  One difference between triassociative algebras and dialgebras is the coboundary operator in cohomology.  Specifically, in both cases, the coboundary involves certain products $\circ^\psi_j$.  The definition of this product is more delicate in the triassociative algebra case than in the dialgebra case.  This is due to the fact that planar trees, which model triassociative algebra cohomology, are more complicated than planar \emph{binary} trees, which model dialgebra cohomology.  A more fundamental difference is that the sets of planar trees form a simplicial set, whereas the sets of binary trees form only an \emph{almost} simplicial set \cite[3.10]{loday}.  This is important, since the author intends to adapt the algebraic works in the present paper to study Topological Triassociative (Co)homology, in analogy with Topological Hochschild (Co)homology (see, e.g., \cite{bok,ekmm,ms}).  Such objects should be defined as the geometric realization of a certain simplicial spectrum (or the $Tot$ of a certain cosimplicial spectrum), modeled after the algebraic (co)simplicial object that defines (co)homology.

In the classical case of an associative algebra $R$, the coefficients in Hochschild homology and cohomology are the same, namely, the $R$-representations.  In contrast, in the case of triassociative algebras (and also dialgebras), the coefficients in homology and cohomology are not the same.  The natural coefficient for triassociative cohomology is a representation, which can be defined using the $11$ triassociative algebra axioms.  This leads to an expected algebraic deformation theory.  To define triassociative homology with coefficients, we will construct the universal enveloping algebra $UA$, which is a unital associative algebra, of a triassociative algebra $A$.  A left $UA$-module is exactly an $A$-representation.  We then define an $A$-\emph{corepresentation} to be a right $UA$-module.  These right $UA$-modules are the coefficients in triassociative homology.  In the classical Hochschild theory, the universal enveloping algebra is $R^e = R \otimes R^{op}$.

We remark that all of the results in this paper can be reproduced for tridendriform and tricubical algebras in \cite{lr} with minimal modifications.

\subsection*{Organization}
The next section is devoted to defining triassociative algebra cohomology with non-trivial coefficients, in which we begin with a brief discussion of representations over a triassociative algebra.  The construction of triassociative algebra cohomology requires a discussion of planar trees.  We observe that the sets of planar trees form a simplicial set (Proposition \ref{prop:simplicial}).  After that, we define the cochain complex $\Ctrias^*(A,M)$ for a triassociative algebra $A$ with coefficients in an $A$-representation $M$ (Theorem \ref{thm:cochain}) and give descriptions of the low dimensional cohomology modules $\Htrias^*$ for $* \leq 2$.  It is also shown that $\Htrias^n(A,-)$ is trivial for all $n \geq 2$ when $A$ is a free triassociative algebra (Theorem \ref{thm:trivial Hi}).

In \S \ref{sec:E(A)}, we construct triassociative algebra homology with non-trivial coefficients by using the \emph{universal enveloping algebra} $UA$ of a triassociative algebra $A$.  The $A$-representations are identified with the left $UA$-modules (Proposition \ref{prop:UA left}).  We then compute the associated graded algebra of $UA$ under the length filtration (Theorem \ref{thm:filtration}).   This leads to a triassociative version of the Poincar\'{e}-Birkhoff-Witt Theorem (Corollary \ref{cor:PBW}).  We then use the \emph{right} $UA$-modules as the coefficients for triassociative algebra homology, analogous to the dialgebra case \cite{frabetti}.  To show that the purported homology complex is actually a chain complex, we relate it to the \emph{cotangent complex} (Proposition \ref{prop:cotangent}), which can be used to define both the homology and the cohomology complexes.  The cotangent complex is the analogue of the bar complex in Hochschild homology.  When the coefficient module is taken to be the ground field (i.e.\ trivial coefficients), our triassociative homology agrees with the one constructed in Loday-Ronco \cite{lr}.  That section ends with descriptions of $\HT_0$ and $\HT_1$.

Section \ref{sec:def} is devoted to studying algebraic deformations of triassociative algebras.  We define algebraic deformations and their infinitesimals for a triassociative algebra $A$.  It is observed that an infinitesimal is always a $2$-cocycle in $\Ctrias^2(A,A)$ whose cohomology class is determined by the equivalence class of the deformation (Theorem \ref{thm:inf}).  A triassociative algebra $A$ is called \emph{rigid} if every deformation of $A$ is equivalent to the trivial deformation.  It is observed that the cohomology module $\Htrias^2(A,A)$ can be thought of as the obstruction to the rigidity of $A$.  Namely, we observe that $A$ is rigid, provided that the module $\Htrias^2(A,A)$ is trivial (Corollary \ref{cor:rigid}).  As examples, free triassociative algebras are rigid (Corollary \ref{cor:free}).  Finally, we identify the obstructions to extending $2$-cocycles in $\Ctrias^2(A,A)$ to deformations.  Given a $2$-cocycle, there is a sequence of obstruction classes in $\Ctrias^3(A,A)$, which are shown to be $3$-cocycles (Lemma \ref{lem:ob}).  The simultaneous vanishing of their cohomology classes is equivalent to the existence of a deformation whose infinitesimal is the given $2$-cocycle (Theorem \ref{thm:ob}).  In particular, these obstructions always vanish if the cohomology module $\Htrias^3(A,A)$ is trivial (Corollary \ref{cor:ob}).  We remark that the work of Balavoine \cite{bal} gives another approach to studying deformations of triassociative algebras.


\section{Cohomology of triassociative algebras}
\label{sec:def complex}

For the rest of this paper, we work over a fixed ground field $K$.  Tensor products are taken over $K$.  We begin by recalling some relevant definitions about triassociative algebras and planar trees from \cite{lr}.

\subsection{Triassociative algebras and representations}
\label{subsec:trias}

A \emph{triassociative algebra} is a vector space $A$ that comes equipped with three binary operations, $\lprod$ (left), $\rprod$ (right), and $\mprod$ (middle), satisfying the following $11$ triassociative axioms for all $x, y, z \in A$:
   \begin{subequations}
   \allowdisplaybreaks
   \label{eq:axioms}
   \begin{align}
   (x \lprod y) \lprod z & = x \lprod (y \lprod z), \label{axiom1} \\
   (x \lprod y) \lprod z & = x \lprod (y \rprod z), \label{axiom2} \\
   (x \rprod y) \lprod z & = x \rprod (y \lprod z), \label{axiom3} \\
   (x \lprod y) \rprod z & = x \rprod (y \rprod z), \label{axiom4} \\
   (x \rprod y) \rprod z & = x \rprod (y \rprod z), \label{axiom5} \\
   (x \lprod y) \lprod z & = x \lprod (y \,\mprod\, z), \label{axiom6} \\
   (x \,\mprod\, y) \lprod z & = x \,\mprod\, (y \lprod z), \label{axiom7} \\
   (x \lprod y) \,\mprod\, z & = x \,\mprod\, (y \rprod z), \label{axiom8} \\
   (x \rprod y) \,\mprod\, z & = x \rprod (y \,\mprod\, z), \label{axiom9} \\
   (x \,\mprod\, y) \rprod z & = x \rprod (y \rprod z), \label{axiom10} \\
   (x \,\mprod\, y) \,\mprod\, z & = x \,\mprod\, (y \,\mprod\, z). \label{axiom11}
   \end{align}
   \end{subequations}
Note that the first $5$ axioms state that $(A, \lprod, \rprod)$ is a dialgebra \cite{loday}, and they do not involve the middle product $\mprod$.  From now on, $A$ will always denote an arbitrary triassociative algebra, unless stated otherwise.

A \emph{morphism} of triassociative algebras is a vector space map that respects the three products.

An \emph{$A$-representation} is a vector space $M$ together with: (i) $3$ left operations $\lprod, \, \rprod, \, \mprod \colon A \otimes M \to M$ and (ii) $3$ right operations $\lprod, \, \rprod, \, \mprod \colon M \otimes A \to M$ satisfying \eqref{axiom1} -- \eqref{axiom11} whenever exactly one of $x, y, z$ is from $M$ and the other two from $A$.  Thus, there are $33$ total axioms.  From now on, $M$ will always denote an arbitrary $A$-representation, unless stated otherwise.

For example, if $\varphi \colon A \to B$ is a morphism of triassociative algebras, then $B$ becomes an $A$-representation via $\varphi$, namely, $x \star b = \varphi(x) \star b$ and $b \star x = b \star \varphi(x)$ for $x \in A$, $b \in B$, and $\star \in \lbrace \lprod, \rprod, \mprod \rbrace$.  In particular, $A$ is naturally an $A$-representation via the identity map.

In order to construct triassociative algebra cohomology, we also need to use planar trees.

\subsection{Planar trees}
\label{subsec:trees}

For integers $n \geq 0$, let $T_n$ denote the set of planar trees with $n + 1$ leaves and one root in which each internal vertex has valence at least $2$.  We will call them \emph{trees} from now on.  The first four sets $T_n$ are listed below:
   \[
   T_0 = \lbrace \begin{picture}(2,10)      
                 \drawline(1,0)(1,10)
                 \end{picture}\,
         \rbrace, \qquad
   T_1 = \lbrace \begin{picture}(16,10)     
                 \drawline(8,0)(8,2)(0,10)
                 \drawline(8,2)(16,10)
                 \end{picture}
         \rbrace, \qquad
   T_2 = \lbrace \begin{picture}(16,10)     
                 \drawline(8,0)(8,2)(0,10)
                 \drawline(8,2)(16,10)
                 \drawline(12,6)(8,10)
                 \end{picture},\,
                 \begin{picture}(16,10)     
                 \drawline(8,0)(8,2)(0,10)
                 \drawline(4,6)(8,10)
                 \drawline(8,2)(16,10)
                 \end{picture},\,
                 \begin{picture}(16,10)     
                 \drawline(8,0)(8,10)
                 \drawline(8,2)(0,10)
                 \drawline(8,2)(16,10)
                 \end{picture}
         \rbrace,
   \]
   \[
   T_3 = \lbrace \begin{picture}(18,11)    
                 \drawline(9,0)(9,2)(0,11)
                 \drawline(9,2)(18,11)
                 \drawline(12,5)(6,11)
                 \drawline(15,8)(12,11)
                 \end{picture},\,
                 \begin{picture}(18,11)    
                 \drawline(9,0)(9,2)(0,11)
                 \drawline(9,2)(18,11)
                 \drawline(12,5)(6,11)
                 \drawline(9,8)(12,11)
                 \end{picture},\,
                 \begin{picture}(18,11)    
                 \drawline(9,0)(9,2)(0,11)
                 \drawline(3,8)(6,11)
                 \drawline(9,2)(18,11)
                 \drawline(15,8)(12,11)
                 \end{picture},\,
                 \begin{picture}(18,11)    
                 \drawline(9,0)(9,2)(0,11)
                 \drawline(6,5)(12,11)
                 \drawline(9,8)(6,11)
                 \drawline(9,2)(18,11)
                 \end{picture},\,
                 \begin{picture}(18,11)    
                 \drawline(9,0)(9,2)(0,11)
                 \drawline(3,8)(6,11)
                 \drawline(6,5)(12,11)
                 \drawline(9,2)(18,11)
                 \end{picture},\,
                 \begin{picture}(18,11)    
                 \drawline(9,0)(9,2)(0,11)
                 \drawline(9,2)(18,11)
                 \drawline(12,5)(6,11)
                 \drawline(12,5)(12,11)
                 \end{picture},\,
                 \begin{picture}(18,11)    
                 \drawline(9,0)(9,2)(0,11)
                 \drawline(9,2)(18,11)
                 \drawline(9,2)(9,11)
                 \drawline(15,8)(12,11)
                 \end{picture},\,
                 \begin{picture}(18,11)    
                 \drawline(9,0)(9,2)(0,11)
                 \drawline(9,2)(18,11)
                 \drawline(9,2)(9,8)
                 \drawline(9,8)(6,11)
                 \drawline(9,8)(12,11)
                 \end{picture},\,
                 \begin{picture}(18,11)    
                 \drawline(9,0)(9,2)(0,11)
                 \drawline(9,2)(18,11)
                 \drawline(9,2)(9,11)
                 \drawline(3,8)(6,11)
                 \end{picture},\,
                 \begin{picture}(18,11)    
                 \drawline(9,0)(9,2)(0,11)
                 \drawline(9,2)(18,11)
                 \drawline(6,5)(6,11)
                 \drawline(6,5)(12,11)
                 \end{picture},\,
                 \begin{picture}(18,11)    
                 \drawline(9,0)(9,2)(0,11)
                 \drawline(9,2)(18,11)
                 \drawline(9,2)(6,11)
                 \drawline(9,2)(12,11)
                 \end{picture}
         \rbrace
   \]
Trees in $T_n$ are said to have \emph{degree $n$}, denoted by $\vert \psi \vert = n$ for $\psi \in T_n$.   The $n + 1$ leaves of a tree in $T_n$ are labeled $\lbrace 0, 1, \ldots , n \rbrace$ from left to right.

A leaf is said to be \emph{left oriented} (respectively, \emph{right oriented}) if it is the left most (respectively, right most) leaf of the vertex underneath it.  Leaves that are neither left nor right oriented are called \emph{middle leaves}.  For example, in the tree $\begin{picture}(18,11) \drawline(9,0)(9,2)(0,11) \drawline(9,2)(18,11) \drawline(9,2)(9,11) \drawline(15,8)(12,11) \end{picture}$, leaves $0$ and $2$ are left oriented, while leaf $3$ is right oriented.  Leaf $1$ is a middle leaf.

Given trees $\psi_0, \ldots , \psi_k$, one can form a new tree by the operation of \emph{grafting}.  Namely, the grafting of these $k + 1$ trees is the tree $\psi_0 \vee \cdots \vee \psi_k$ obtained by arranging $\psi_0, \ldots , \psi_k$ from left to right and joining the $k + 1$ roots to form a new (lowest) internal vertex, which is connected to a new root.  The degree of $\psi_0 \vee \cdots \vee \psi_k$ is
   \[
   \lbrace (\vert \psi_0 \vert + 1) + \cdots + (\vert \psi_k \vert + 1) \rbrace - 1
   \,=\, k + \sum_{i=0}^k \vert \psi_i \vert.
   \]
Conversely, every tree $\psi$ can be written uniquely as the grafting of $k + 1$ trees, $\psi_0 \vee \cdots \vee \psi_k$, where the valence of the lowest internal vertex of $\psi$ is $k + 1$.

Before defining cohomology, we make the following observation, which is not used in the rest of the paper but maybe useful in future studies of triassociative cohomology.

For $0 \leq i \leq n+1$, define a function
   \[
   d_i \colon T_{n+1} \to T_n,
   \]
which sends a tree $\psi \in T_{n+1}$ to the tree $d_i \psi \in T_n$ obtained from $\psi$ by deleting the $i$th leaf.  For $0 \leq i \leq n$, define another function
   \[
   s_i \colon T_n \to T_{n+1}
   \]
as follows: For $\psi \in T_n$, $s_i \psi \in T_{n+1}$ is the tree obtained from $\psi$ by adding a new leaf to the internal vertex connecting to leaf $i$, and this new leaf is placed immediately to the left of the original leaf $i$.  For example, if $\psi = \begin{picture}(16,10)     
                 \drawline(8,0)(8,2)(0,10)
                 \drawline(8,2)(16,10)
                 \drawline(12,6)(8,10)
                 \end{picture}$,
then $s_0(\psi) = \begin{picture}(18,11)    
                 \drawline(9,0)(9,2)(0,11)
                 \drawline(9,2)(18,11)
                 \drawline(9,2)(9,11)
                 \drawline(15,8)(12,11)
                 \end{picture}$
and $s_1(\psi) = s_2(\psi) = \begin{picture}(18,11)    
                 \drawline(9,0)(9,2)(0,11)
                 \drawline(9,2)(18,11)
                 \drawline(12,5)(6,11)
                 \drawline(12,5)(12,11)
                 \end{picture}$.

\begin{prop}
\label{prop:simplicial}
The sets $\lbrace T_n \rbrace_{n \geq 0}$ form a simplicial set with face maps $d_i$ and degeneracy maps $s_i$.
\end{prop}

\begin{proof}
Recall that the simplicial relations are:
   \[
   \begin{split}
   d_id_j & \,=\, d_{j-1}d_i, \quad \text{if }i < j, \\
   d_is_j & \,=\, \begin{cases} s_{j-1}d_i & \text{if } i < j, \\
                                \Id &\text{if } i = j, j+1 \\
                                s_jd_{i-1} &\text{if } i > j+1, \end{cases} \\
   s_is_j & \,=\, s_{j+1}s_i, \quad \text{if }i \leq j.
   \end{split}
   \]
All of them are immediate from the definitions.
\end{proof}

\subsection{Cohomology}
\label{subsec:cohom}

For integers $n \geq 0$, define the \emph{module of $n$-cochains of $A$ with coefficients in $M$} to be
   \[
   \Ctrias^n(A,M) \,:=\, \Hom_K(K \lbrack T_n \rbrack \otimes A^{\otimes n}, M).
   \]
To define the coboundary maps, we need the following operations.  For $0 \leq i \leq n+1$, define a function
   \[
   \circ_i \colon T_{n+1} \to \lbrace \lprod, \rprod, \mprod \rbrace
   \]
according to the following rules.  Let $\psi$ be a tree in $T_{n+1}$, which is written uniquely as $\psi = \psi_0 \vee \cdots \vee \psi_k$ in which the valence of the lowest internal vertex of $\psi$ is $k + 1$.  Also, write $\circ_i^\psi$ for $\circ_i(\psi)$.  Then $\circ^\psi_0$ is given by
   \[
   \circ^\psi_0 \,=\,
   \begin{cases}
   \lprod & \text{if } \vert \psi_0 \vert = 0 \text{ and } k = 1, \\
   \rprod & \text{if } \vert \psi_0 \vert > 0, \\
   \mprod & \text{if } \vert \psi_0 \vert = 0 \text{ and } k > 1.
   \end{cases}
   \]
For $1 \leq i \leq n$, $\circ^\psi_i$ is given by
   \[
   \circ^\psi_i \,=\,
   \begin{cases}
   \lprod & \text{if the } i\text{th leaf of } \psi \text{ is left oriented}, \\
   \rprod & \text{if the } i\text{th leaf of } \psi \text{ is right oriented}, \\
   \mprod & \text{if the } i\text{th leaf of } \psi \text{ is a middle leaf}.
   \end{cases}
   \]
Finally, $\circ^\psi_{n+1}$ is given by
   \[
   \circ^\psi_{n+1} \,=\,
   \begin{cases}
   \lprod & \text{if } \vert \psi_k \vert > 0, \\
   \rprod & \text{if } k = 1 \text{ and } \vert \psi_1 \vert = 0, \\
   \mprod & \text{if } k > 1 \text{ and } \vert \psi_k \vert = 0.
   \end{cases}
   \]
For example, for the $11$ trees $\psi$ in $T_3$ (from left to right), we have:
   \[
   \begin{cases}
   \circ^\psi_0 &  = ~ \lprod, \lprod, \rprod, \rprod, \rprod, \lprod, \mprod, \mprod, \rprod, \rprod, \mprod, \\
   \circ^\psi_{n+1} &  = ~ \lprod, \lprod, \lprod, \rprod, \rprod, \lprod, \lprod, \mprod, \mprod, \rprod, \mprod.
   \end{cases}
   \]
Now define the map
   \[
   \delta^n \colon \Ctrias^n(A,M) \,\to\, \Ctrias^{n+1}(A,M)
   \]
to be the alternating sum,
   \[
   \delta^n \,=\, \sum_{i=0}^{n+1} \, (-1)^i \delta^n_i,
   \]
where
   \[
   (\delta^n_i f)(\psi; a_1, \ldots , a_{n+1}) \,=\,
   \begin{cases}
   a_1 \circ^\psi_0 f(d_0 \psi; a_2, \ldots , a_{n+1}) & \text{if } i = 0, \\
   f(d_i \psi; a_1, \ldots , a_i \circ^\psi_i a_{i+1}, \ldots , a_n) & \text{if } 1 \leq i \leq n, \\
   f(d_{n+1}\psi ; a_1, \ldots , a_n) \circ^\psi_{n+1} a_{n+1} & \text{if } i = n+1,
   \end{cases}
   \]
for $f \in \Ctrias^n(A,M)$, $\psi \in T_{n+1}$, and $a_1, \ldots , a_{n+1} \in A$.

\begin{thm}
\label{thm:cochain}
The maps $\delta^n_i$ satisfy the cosimplicial identities,
   \[
   \delta^{n+1}_j \delta^n_i \,=\, \delta^{n+1}_i \delta^n_{j-1},
   \]
for $0 \leq i < j \leq n + 2$.  In particular, $(\Ctrias^*(A,M), \delta)$ is a cochain complex.
\end{thm}

\begin{proof}
First note that for $1 \leq i \leq n$, the map $\delta^n_i$ is dual to the map
   \[
   d_i \colon C^{\text{Trias}}_{n+1}(A) = K \lbrack T_{n+1} \rbrack \otimes A^{\otimes n+1} \,\to\, K \lbrack T_n \rbrack \otimes A^{\otimes n} = C^{\text{Trias}}_n(A)
   \]
in Loday-Ronco \cite{lr}, i.e.\ $\delta^n_i = \Hom(d_i, M)$.   It is shown there that these $d_i$ satisfy the simplicial identities.  Therefore, to prove the cosimplicial identities, it suffices to consider the following three cases: (i) $i = 0$, $j = n + 2$, (ii) $i = 0 < j < n + 2$, and (iii) $0 < i < j = n + 2$.

Consider case (i).  Suppose that $\psi \in T_{n+2}$, $a_1, \ldots , a_{n+2} \in A$, and $f \in \Ctrias^n(A,M)$.  Let $\alpha$ be the element $(\psi ; a_1, \ldots , a_{n+2})$ in $K \lbrack T_{n+2} \rbrack \otimes A^{\otimes n + 2}$.  Write $\psi$ as $\psi_0 \vee \cdots \vee \psi_k$ uniquely with $k + 1$ being the valence of the lowest internal vertex of $\psi$.  We have
   \[
   \allowdisplaybreaks
   \begin{split}
   (\delta^{n+1}_{n+2} \, \delta^n_0 f)(\alpha)
   & \,=\, \lbrace (\delta^n_0 f)(d_{n+2}\psi ; a_1, \ldots , a_{n+1}) \rbrace \circ^\psi_{n+2} a_{n+2} \\
   & \,=\, \lbrace a_1 \circ^{d_{n+2}\psi}_0 f(d_0d_{n+2} \psi ; a_2, \ldots , a_{n+1}) \rbrace \circ^\psi_{n+2} a_{n+2} \\
   & \,=\, (a_1 \circ^{d_{n+2}\psi}_0 y) \circ^\psi_{n+2} a_{n+2},
   \end{split}
   \]
where $y = f(d_{n+1} d_0 \psi ; a_2, \ldots , a_{n+1})$.  Here we used the simplicial identity $d_0d_{n+2} = d_{n+1} d_0$.  On the other hand, we have
   \[
   \allowdisplaybreaks
   \begin{split}
   (\delta^{n+1}_0 \delta^n_{n+1} f)(\alpha)
   & \,=\, a_1 \circ^\psi_0 \lbrace (\delta^n_{n+1} f)(d_0 \psi ; a_2, \ldots , a_{n+2}) \rbrace \\
   & \,=\, a_1 \circ^\psi_0 (y \circ^{d_0\psi}_{n+1} a_{n+2}).
   \end{split}
   \]
In order to show that
   \begin{equation}
   \label{eq:i=0}
   (a_1 \circ^{d_{n+2}\psi}_0 y) \circ^\psi_{n+2} a_{n+2}
   \,=\, a_1 \circ^\psi_0 (y \circ^{d_0\psi}_{n+1} a_{n+2}),
   \end{equation}
we need to consider $7$ different cases.

(1) If $k = 1$ and $\vert \psi_0 \vert = 0$, then $\circ^\psi_0 = ~ \lprod ~ = \circ^{d_{n+2}\psi}_0$.  In this case, it follows from the triassociative algebra axioms \eqref{axiom1}, \eqref{axiom2}, and \eqref{axiom6} that
   \[
   (a_1 \lprod y) \ast a_{n+2} \,=\, a_1 \lprod (y \ast^{\prime} a_{n+2})
   \]
regardless of what $\ast, \ast^\prime \in \lbrace \lprod, \rprod, \mprod \rbrace$ are.

(2) If $k = 1$ and $\vert \psi_1 \vert = 0$, then $\circ^\psi_{n+2} = ~ \rprod ~ = \circ^{d_0\psi}_{n+1}$.  It then follows from axioms \eqref{axiom4}, \eqref{axiom5}, and \eqref{axiom10} that
   \[
   (a_1 \ast y) \rprod a_{n+2} \,=\, a_1 \ast^\prime (y \rprod a_{n+2})
   \]
regardless of what $\ast, \ast^\prime \in \lbrace \lprod, \rprod, \mprod \rbrace$ are.

(3) For $k \geq 1$, if both $\vert \psi_0 \vert$ and $\vert \psi_k \vert$ are positive, then $\circ^{d_{n+2}\psi}_0 = ~ \rprod ~ = \circ^\psi_0$ and $\circ^\psi_{n+2} = ~ \lprod ~ = \circ^{d_0 \psi}_{n+1}$.  Therefore, \eqref{eq:i=0} holds by axiom \eqref{axiom3}.

(4) If $k = 2$ and $\vert \psi_0 \vert = 0 = \vert \psi_2 \vert$, then $\circ^{d_{n+2}\psi}_0 = ~\lprod$, $\circ^\psi_{n+2} = ~ \mprod$, $\circ^\psi_0 = ~\mprod$, and $\circ^{d_0\psi}_{n+1} = ~\rprod$.  It follows from axiom \eqref{axiom8} that \eqref{eq:i=0} holds.

(5) If $k > 2$ and $\vert \psi_0 \vert = 0 = \vert \psi_k \vert$, then $\circ^{d_{n+2}\psi}_0 = \circ^\psi_{n+2} = \circ^\psi_0 = \circ^{d_0\psi}_{n+1} = ~ \mprod$.  Thus, \eqref{eq:i=0} follows from axiom \eqref{axiom11}.

(6) If $k \geq 2$, $\vert \psi_0 \vert = 0$, and $\vert \psi_k \vert > 0$, then $\circ^{d_{n+2}\psi}_0 = ~\mprod~ = \circ^\psi_0$ and $\circ^\psi_{n+2} = ~\lprod~ = \circ^{d_0\psi}_{n+1}$.  Therefore, it follows from axiom \eqref{axiom7} that \eqref{eq:i=0} holds.

(7) If $k \geq 2$, $\vert \psi_0 \vert > 0$, and $\vert \psi_k \vert = 0$, then $\circ^{d_{n+2}\psi}_0 = ~\rprod~ = \circ^\psi_0$ and $\circ^\psi_{n+2} = ~\mprod~ = \circ^{d_0\psi}_{n+1}$.  Therefore, \eqref{eq:i=0} follows from axiom \eqref{axiom9}.

Therefore, the cosimplicial identity holds in the case $i = 0$, $j = n + 2$.

Next we consider case (ii), where $i = 0$, $j = 1, \ldots , n + 1$.  We have
   \[
   \begin{split}
   (\delta^{n+1}_j \delta^n_0 f)(\alpha)
   & \,=\, (\delta^n_0 f)(d_j \psi ; a_1, \ldots , a_j \circ^\psi_j a_{j+1}, \ldots , a_{n+2}) \\
   & \,=\, \begin{cases}
           (a_1 \circ^\psi_1 a_2) \circ^{d_1\psi}_0 z & \text{if } j = 1 \\
           a_1 \circ^{d_j\psi}_0 x  & \text{if } 2 \leq j \leq n + 1, \end{cases}
   \end{split}
   \]
where
   \[
   \begin{split}
   z & = f(d_0 d_0 \psi ; a_3, \ldots , a_{n+2}), \\
   x & = f(d_{j-1}d_0 \psi ; a_2, \ldots , a_j \circ^\psi_j a_{j+1}, \ldots , a_{n+2}).
   \end{split}
   \]
Here we used the simplicial identity $d_0 d_j = d_{j-1} d_0$ for $0 < j$.  On the other hand, we have
   \[
   \begin{split}
   (\delta^{n+1}_0 \delta^n_{j-1} f)(\alpha)
   & \,=\, a_1 \circ^\psi_0 (\delta^n_{j-1}f)(d_0 \psi ; a_2, \ldots , a_{n+2}) \\
   & \,=\, \begin{cases} a_1 \circ^\psi_0 (a _2 \circ^{d_0\psi}_0 z) & \text{if } j = 1 \\
                     a_1 \circ^\psi_0 w & \text{if } 2 \leq j \leq n + 1, \end{cases}
   \end{split}
   \]
where
   \[
   w \,=\, f(d_{j-1}d_0 \psi ; a_2, \ldots , a_j \circ^{d_0\psi}_{j-1} a_{j+1}, \ldots , a_{n+2}).
   \]
In the cases $2 \leq j \leq n + 1$, we have that $\circ^\psi_j = \circ^{d_0 \psi}_{j-1}$ and $\circ^{d_j\psi}_0 = \circ^\psi_0$, and hence $a_1 \circ^{d_j\psi}_0 x = a_1 \circ^\psi_0 w$.  In the case $j = 1$, we need to show the identity
   \begin{equation}
   \label{eq:cohom j=1}
   (a_1 \circ^\psi_1 a_2) \circ^{d_1\psi}_0 z
   \,=\, a_1 \circ^\psi_0 (a _2 \circ^{d_0\psi}_0 z).
   \end{equation}
We break it into three cases.

(1) If $\vert \psi_0 \vert \geq 2$, then $\circ^{d_1\psi}_0 = ~\rprod~ = \circ^{d_0\psi}_0$.  Therefore, it follows from axioms \eqref{axiom4}, \eqref{axiom5}, \eqref{axiom10} that
   \[
   (a_1 \ast a_2) \rprod z \,=\, a_1 \ast^\prime (a_2 \rprod z)
   \]
regardless of what $\ast, \ast^\prime \in \lbrace \lprod, \rprod, \mprod \rbrace$ are.

(2) If $\vert \psi_0 \vert = 1$, then $\circ^\psi_1 = ~\rprod~ = \circ^\psi_0$ and
   \[
   \circ^{d_0 \psi}_0 \,=\, \circ^{d_1\psi}_0 \,=\,
   \begin{cases} \lprod & \text{if } k = 1, \\
                 \mprod & \text{if } k \geq 2.
   \end{cases}
   \]
Therefore, \eqref{eq:cohom j=1} holds by axiom \eqref{axiom3} if $k = 1$ and by axiom \eqref{axiom9} if $k \geq 2$.

(3) Now suppose that $\vert \psi_0 \vert = 0$.  If $k = 1$, then $\circ^\psi_0 = ~\lprod~ = \circ^\psi_1$.  It follows that \eqref{eq:cohom j=1} holds by axioms \eqref{axiom1}, \eqref{axiom2}, and \eqref{axiom6}.  If $k \geq 2$, then $\circ^\psi_0 = \mprod$.  To figure out what the other three operations are, we need to consider two sub-cases.
   \begin{itemize}
   \item If $\vert \psi_1 \vert > 0$, then $\circ^{d_0\psi}_0 = ~\rprod$, $\circ^\psi_1 = ~\lprod$, and $\circ^{d_1\psi}_0 = \mprod$.  Thus, \eqref{eq:cohom j=1} holds by axiom \eqref{axiom8}.
   \item If $\vert \psi_1 \vert = 0$, then $\circ^\psi_1 = \mprod$ and
     \[
     \circ^{d_1\psi}_0 \,=\, \circ^{d_0\psi}_0 \,=\, \begin{cases} \lprod & \text{if } k=2, \\ \mprod & \text{if } k > 2. \end{cases}
     \]
It follows that \eqref{eq:cohom j=1} holds by axiom \eqref{axiom7} when $k = 2$ and by axiom \eqref{axiom11} when $k > 2$.
   \end{itemize}
This proves the cosimplicial identities for the $\delta^n_l$ when $i = 0$ and $1 \leq j \leq n + 1$.  The proof for the case $1 \leq i \leq n + 1$, $j = n + 2$ is similar to the argument that was just given.
\end{proof}

In view of Theorem \ref{thm:cochain}, we define the \emph{$n$th cohomology of $A$ with coefficients in the representation $M$} to be
   \[
   \Htrias^n(A,M) \,:=\, H^n(\Ctrias^*(A,M),\delta)
   \]
for $n \geq 0$.  We describe the first three cohomology modules below.

\subsection{$\Htrias^0$ and $\Htrias^1$}
\label{subsec:H0 H1}

A linear map $\varphi \colon A \to M$ is called a \emph{derivation of $A$ with values in $M$} if it satisfies the condition,
   \[
   \varphi(a \ast b) \,=\, \varphi(a) \ast b + a \ast \varphi(b),
   \]
for all $a, b \in A$ and $\ast \in \lbrace \lprod, \rprod, \mprod \rbrace$.  Denote by $\Der(A,M)$ the submodule of $\Hom_K(A,M)$ consisting of all the derivations of $A$ with values in $M$.  For an element $m \in M$, define the map $ad_m \colon A \to M$, where
   \[
   ad_m(a) \,:=\, a \lprod m - m \rprod a
   \]
for $a \in A$.  Such a map is called an \emph{inner derivation of $A$ with values in $M$}.  It follows immediately from the triassociative algebra axioms \eqref{eq:axioms} that each map $ad_m$ belongs to $\Der(A,M)$.   Let $\Inn(A,M)$ denote the submodule of $\Der(A,M)$ consisting of all the inner derivations of $A$ with values in $M$.

Identify $\Ctrias^0(A,M)$ with $M$ and $\Ctrias^1(A,M)$ with $\Hom_K(A,M)$.  Under such identifications, the coboundary map $\delta^0 \colon M \to \Hom_K(A,M)$ is given by $\delta^0(m) \,=\, ad_m.$  Therefore, the cohomology module $\Htrias^0(A,M)$, which is the kernel of $\delta^0$, is the following submodule of $M$:
   \[
   \Htrias^0(A,M) \,\cong\, M^A \,:=\, \lbrace m \in M \, \vert \, a \lprod m = m \rprod a \text{ for all } a \in A \rbrace.
   \]

The image of $\delta^0$ is the module $\Inn(A,M)$.  Now if $f \in \Ctrias^1(A,M)$, then
   \[
   (\delta^1 f)(\psi ; a, b) \,=\, a \ast f(b) - f(a \ast b) + f(a) \ast b,
   \]
where
   \begin{equation}
   \label{eq:product}
   \ast \,=\,
   \begin{cases}
   \lprod & \text{if } \psi = \begin{picture}(16,10)     
                 \drawline(8,0)(8,2)(0,10)
                 \drawline(8,2)(16,10)
                 \drawline(12,6)(8,10)
                 \end{picture}, \\
   \rprod & \text{if } \psi = \begin{picture}(16,10)     
                 \drawline(8,0)(8,2)(0,10)
                 \drawline(4,6)(8,10)
                 \drawline(8,2)(16,10)
                 \end{picture}, \\
   \mprod & \text{if } \psi = \begin{picture}(16,10)     
                 \drawline(8,0)(8,10)
                 \drawline(8,2)(0,10)
                 \drawline(8,2)(16,10)
                 \end{picture}.

   \end{cases}
   \end{equation}
In particular, the kernel of $\delta^1$ is exactly $\Der(A,M)$.  Therefore, as in the cases of associative algebras and dialgebras \cite{frabetti}, we have
   \[
   \Htrias^1(A,M) \,\cong\, \frac{\Der(A,M)}{\Inn(A,M)}.
   \]

\subsection{$\Htrias^i(A,-)$ when $A$ is free}
\label{subsec:Hi}

It is shown in \cite[Theorem 4.1]{lr} that the operad $Trias$ for triassociative algebras is Koszul.  In other words, $H^{\text{Trias}}_n(A,K) = 0$ for $n \geq 2$ if $A$ is a free triassociative algebra.  Here $H^{\text{Trias}}_n(A,K)$ denotes the triassociative homology of $A$ with trivial coefficients defined by Loday and Ronco \cite{lr}.  It agrees with our triassociative homology defined below.  The Koszulness of $Trias$ implies that free triassociative algebras have trivial higher cohomology modules.

\begin{thm}
\label{thm:trivial Hi}
Let $A$ be a free triassociative algebra.  Then 
   \[
   \Htrias^n(A,M) = 0
   \]
for all $n \geq 2$ and any $A$-representation $M$.
\end{thm}

\begin{proof}
The first line in the proof of Theorem \ref{thm:cochain} implies that, for any triassociative algebra $A$,
   \[
   \Htrias^n(A,M) \cong \Hom_K(H^{\text{Trias}}_n(A,K),M).
   \]
The Theorem now follows from the Koszulness of $Trias$.
\end{proof}

\subsection{Abelian extensions and $\Htrias^2$}
\label{subsec:H2}

Define an \emph{abelian triassociative algebra} to be a triassociative algebra $P$ in which all three products, $\lprod$, $\rprod$, $\mprod$, are equal to $0$.  In this case, we will just say that $P$ is abelian.  Any vector space becomes an abelian triassociative algebra when equipped with the trivial products.   Suppose that $\xi \colon 0 \to P \xrightarrow{i} E \xrightarrow{\pi} A \to 0$ is a short exact sequence of triassociative algebras in which $P$ is abelian.  Then $P$ has an induced $A$-representation structure via $a \ast p  = e \ast i(p)$ and $p \ast a  = i(p) \ast e$ for $\ast \in \lbrace \lprod, \rprod, \mprod \rbrace$, $p \in P$, $a \in A$, and any element $e \in E$ such that $\pi(e) = a$.

Now consider an $A$-representation $M$.  By an \emph{abelian extension of $A$ by $M$}, we mean a short exact sequence
   \begin{equation}
   \label{eq:abelian ext}
   \xi\colon \quad 0 \xrightarrow{} M \xrightarrow{i} E \xrightarrow{\pi} A \xrightarrow{} 0
   \end{equation}
of triassociative algebras in which $M$ is abelian and such that the induced $A$-representation structure on $M$ coincides with the original one.  An abelian extension is said to be \emph{trivial} if it splits triassociative algebras.  Given another abelian extension $\xi^\prime = \lbrace 0 \to M \to E^\prime \to A \to 0 \rbrace$ of $A$ by $M$, we say that $\xi$ and $\xi^\prime$ are \emph{equivalent} if there exists a map $\varphi \colon E \to E^\prime$ of triassociative algebras making the obvious ladder diagram commutative.  (Note that such a map $\varphi$ must be an isomorphism.)  Denote by $\lbrack \xi \rbrack$ the equivalence class of an abelian extension $\xi$ and by $\Ext(A,M)$ the set of equivalence classes of abelian extensions of $A$ by $M$.

Suppose that $\xi$ is an abelian extension of $A$ by $M$ as in \eqref{eq:abelian ext}.  By choosing a vector space splitting $\sigma \colon A \to E$, one can identify the underlying vector space of $E$ with $M \oplus A$.  As usual, there exists a map
   \[
   f_\xi \colon K \lbrack T_2 \rbrack \otimes A^{\otimes 2} \,\to\, M
   \]
such that the products in $E$ become
   \begin{equation}
   \label{eq:product E}
   (m,a) \ast (n,b) \,=\,
   (m \ast b + a \ast n + f_\xi(\psi; a, b),\, a \ast b)
   \end{equation}
for $\ast \in \lbrace \lprod, \rprod, \mprod \rbrace$, $m, n \in M$, and $a, b \in A$.  Here $\psi \in T_2$ is given by
   \begin{equation}
   \label{eq:E}
   \psi \,=\,
   \begin{cases}
   \begin{picture}(16,10)     
   \drawline(8,0)(8,2)(0,10)
   \drawline(8,2)(16,10)
   \drawline(12,6)(8,10)
   \end{picture} & \text{if } \ast = ~\lprod, \\
   \begin{picture}(16,10)     
   \drawline(8,0)(8,2)(0,10)
   \drawline(4,6)(8,10)
   \drawline(8,2)(16,10)
   \end{picture} & \text{if } \ast = ~ \rprod, \\
   \begin{picture}(16,10)     
   \drawline(8,0)(8,10)
   \drawline(8,2)(0,10)
   \drawline(8,2)(16,10)
   \end{picture} & \text{if } \ast = \mprod.
   \end{cases}
   \end{equation}
Note that this is basically \eqref{eq:product}.  We will always identify the $3$ trees in $T_2$ with the products $\lbrace \lprod, \rprod, \mprod \rbrace$ like this.  It is easy to check that the $11$ triassociative algebra axioms \eqref{eq:axioms} in $E$ are equivalent to $f_\xi \in \Ctrias^2(A,M)$ being a $2$-cocycle.  For instance, we have
   \[
   \begin{split}
   ((0,x) \lprod (0,y)) \lprod (0,z)
   & = (f_\xi(\begin{picture}(16,10)     
   \drawline(8,0)(8,2)(0,10)
   \drawline(8,2)(16,10)
   \drawline(12,6)(8,10)
   \end{picture} ; x, y), \, x \lprod y) \lprod (0,z) \\
   & = (f_x(\begin{picture}(16,10)     
   \drawline(8,0)(8,2)(0,10)
   \drawline(8,2)(16,10)
   \drawline(12,6)(8,10)
   \end{picture} ; x, y) \lprod z \,+\, f_\xi(\begin{picture}(16,10)     
   \drawline(8,0)(8,2)(0,10)
   \drawline(8,2)(16,10)
   \drawline(12,6)(8,10)
   \end{picture} ; x \lprod y, z), \, (x \lprod y) \lprod z).
   \end{split}
   \]
On the other hand, a similar calculation yields
   \[
   \begin{split}
   (0,x) \lprod ((0,y) \,\mprod\, (0,z))
   & = (x \lprod f_\xi(\begin{picture}(16,10)     
   \drawline(8,0)(8,10)
   \drawline(8,2)(0,10)
   \drawline(8,2)(16,10)
   \end{picture} ; y, z) + f_\xi(\begin{picture}(16,10)     
   \drawline(8,0)(8,2)(0,10)
   \drawline(8,2)(16,10)
   \drawline(12,6)(8,10)
   \end{picture} ; x, y \,\mprod\, z), \, x \lprod (y \,\mprod\, z)).
   \end{split}
   \]
These two expressions are equal by axiom \eqref{axiom6}.  By equating the first factors, we obtain
   \[
   (\delta^2 f_\xi)(\begin{picture}(18,11)    
                 \drawline(9,0)(9,2)(0,11)
                 \drawline(9,2)(18,11)
                 \drawline(12,5)(6,11)
                 \drawline(12,5)(12,11)
                 \end{picture} ; x, y, z) \,=\,0.
   \]
Similar arguments, using the other $10$ triassociative algebra axioms, show that $(\delta^2 f_\xi)(\psi ; x, y, z) = 0$ for the other $10$ trees $\psi \in T_3$.

Conversely, suppose that $g \in \Ctrias^2(A,M)$ is a $2$-cocycle.  Then one can define a triassociative algebra structure on the vector space $M \oplus A$ using \eqref{eq:product E} with $g$ in place of $f_\xi$.  Again, the triassociative algebra axioms are verified because of the cocycle condition on $g$.  This yields an abelian extension of $A$ by $M$,
   \[
   \zeta_g \colon \quad
   0 \xrightarrow{} M \xrightarrow{i} M \oplus A \xrightarrow{\pi} A \xrightarrow{} 0
   \]
in which $i$ and $\pi$ are, respectively, the inclusion into the first factor and the projection onto the second factor.

\begin{thm}
\label{thm:H2}
The above constructions induce well-defined maps
   \begin{align*}
   \Ext(A,M) & \to \Htrias^2(A,M) & \lbrack \xi \rbrack  \mapsto \lbrack f_\xi \rbrack \\
   \Htrias^2(A,M) & \to \Ext(A,M) & \lbrack g \rbrack \mapsto \lbrack \zeta_g \rbrack
   \end{align*}
that are inverse to each other.  In particular, there is a canonical bijection $\Htrias^2(A,M) \cong \Ext(A,M)$.
\end{thm}

The proof is basically identical to that of the classical case for associative algebras (see, e.g., \cite[p.311-312]{weibel} or \cite[2.8-2.9]{frabetti}).  Therefore, we omit the details.  Note that one only needs to see that the maps are well-defined, since they are clearly inverse to each other.


\section{Universal enveloping algebra and triassociative homology}
\label{sec:E(A)}

The purpose of this section is to construct triassociative algebra homology with non-trivial coefficients.  In \cite{lr} Loday and Ronco already constructed triassociative homology with trivial coefficients (i.e.\ $K$).  Our homology agrees with the one in Loday-Ronco \cite{lr} by taking coefficients in $K$.  In order to obtain the non-trivial coefficients, we first need to discuss the universal enveloping algebra.  The dialgebra analogue of the results in this section is worked out in \cite{frabetti}.

\subsection{Universal enveloping algebra}
\label{subsec:uni}

Fix a triassociative algebra $A$.  We would like to identify the $A$-representations as the left modules of a certain associative algebra.  This requires $6$ copies of $A$, since $3$ copies are needed for the left actions and another $3$ for the right actions.  Therefore, we make the following definition.  For a $K$-vector space $V$, $T(V)$ denotes the tensor algebra of $V$, which is the free unital associative $K$-algebra generated by $V$.

Define the \emph{universal enveloping algebra} of $A$ to be the unital associative $K$-algebra $UA$ obtained from the tensor algebra $T(\alpha_l A \oplus \alpha_r A \oplus \alpha_m A \oplus \beta_l A \oplus \beta_r A \oplus \beta_m A)$ on $6$ copies of $A$ by imposing the following $33$ relations for $a, b \in A$ ($3$ relations for each of the $11$ triassociative axioms \eqref{eq:axioms}):
   \[
   \allowdisplaybreaks
   \begin{split}
   &
   \begin{cases}
   \beta_l(b) \cdot \beta_l(a) \buildrel (1) \over = \beta_l(a \lprod b) \buildrel (2) \over = \beta_l(a \,\mprod\, b) \buildrel (3) \over = \beta_l(a \rprod b) & \\
   \beta_l(b) \cdot \alpha_l(a) \buildrel (4) \over = \alpha_l(a) \cdot \beta_l(b) \buildrel (5) \over = \alpha_l(a) \cdot \beta_m(b) \buildrel (6) \over = \alpha_l(a) \cdot \beta_r(b) & \\
   \alpha_l(a \lprod b) \buildrel (7) \over = \alpha_l(a) \cdot \alpha_l(b) \buildrel (8) \over = \alpha_l(a) \cdot \alpha_m(b) \buildrel (9) \over = \alpha_l(a) \cdot \alpha_r(b)
   \end{cases} \\
   &
   \begin{cases}
   \beta_r(a \rprod b) \buildrel (10) \over = \beta_r(b) \cdot \beta_r(a) \buildrel (11) \over = \beta_r(b) \cdot \beta_l(a) \buildrel (12) \over = \beta_r(b) \cdot \beta_m(a) & \\
   \alpha_r(a) \cdot \beta_r(b) \buildrel (13) \over = \beta_r(b) \cdot \alpha_r(a) \buildrel (14) \over = \beta_r(b) \cdot \alpha_l(a) \buildrel (15) \over = \beta_r(b) \cdot \alpha_m(a) & \\
   \alpha_r(a) \cdot \alpha_r(b) \buildrel (16) \over = \alpha_r(a \rprod b) \buildrel (17) \over = \alpha_r(a \lprod b) \buildrel (18) \over = \alpha_r(a \,\mprod\, b)
   \end{cases} \\
   &
   \begin{cases}
   \beta_l(b) \cdot \beta_r(a) \buildrel (19) \over = \beta_r(a \lprod b) & \\
   \beta_l(b) \cdot \alpha_r(a) \buildrel (20) \over = \alpha_r(a) \cdot \beta_l(b) & \\
   \alpha_l(a \rprod b) \buildrel (21) \over = \alpha_r(a) \cdot \alpha_l(b)
   \end{cases}
   \begin{cases}
   \beta_l(b) \cdot \beta_m(a) \buildrel (22) \over = \beta_m(a \lprod b) & \\
   \beta_l(b) \cdot \alpha_m(a) \buildrel (23) \over = \alpha_m(a) \cdot \beta_l(b) & \\
   \alpha_l(a \,\mprod\, b) \buildrel (24) \over = \alpha_m(a) \cdot \alpha_l(b)
   \end{cases}
   \end{split}
   \]

   \[
   \begin{split}
   &
   \begin{cases}
   \beta_m(b) \cdot \beta_l(a) \buildrel (25) \over = \beta_m(a \rprod b) & \\
   \beta_m(b) \cdot \alpha_l(a) \buildrel (26) \over = \alpha_m(a) \cdot \beta_r(b) & \\
   \alpha_m(a \lprod b) \buildrel (27) \over = \alpha_m(a) \cdot \alpha_r(b)
   \end{cases}
   \begin{cases}
   \beta_m(b) \cdot \beta_r(a) \buildrel (28) \over = \beta_r(a \,\mprod\, b) & \\
   \beta_m(b) \cdot \alpha_r(a) \buildrel (29) \over = \alpha_r(a) \cdot \beta_m(b) & \\
   \alpha_m(a \rprod b) \buildrel (30) \over = \alpha_r(a) \cdot \alpha_m(b)
   \end{cases} \\
   &
   \begin{cases}
   \beta_m(b) \cdot \beta_m(a) \buildrel (31) \over = \beta_m(a \,\mprod\, b) & \\
   \beta_m(b) \cdot \alpha_m(a) \buildrel (32) \over = \alpha_m(a) \cdot \beta_m(b) & \\
   \alpha_m(a \,\mprod\, b) \buildrel (33) \over = \alpha_m(a) \cdot \alpha_m(b)
   \end{cases}
   \end{split}
   \]

\begin{prop}
\label{prop:UA left}
Let $M$ be a $K$-vector space.  Then an $A$-representation structure on $M$ is equivalent to a left $UA$-module structure on $M$.
\end{prop}

\begin{proof}
The correspondence is given by ($a \in A$, $x \in M$):
   \[
   \begin{cases}
   \alpha_l(a) \cdot x = a \lprod x & \\
   \alpha_r(a) \cdot x = a \rprod x & \\
   \alpha_m(a) \cdot x = a \,\mprod\, x,
   \end{cases}
   \begin{cases}
   \beta_l(a) \cdot x = x \lprod a & \\
   \beta_r(a) \cdot x = x \rprod a & \\
   \beta_m(a) \cdot x = x \,\mprod\, a.
   \end{cases}
   \]
So $\alpha_*$ and $\beta_*$ correspond to left and right $A$-actions, respectively, and the subscript indicates where the product points to.  In particular, the conditions (1) - (9) above correspond to axioms \eqref{axiom1}, \eqref{axiom2}, and \eqref{axiom6}.  Conditions (10) - (18) correspond to axioms \eqref{axiom4}, \eqref{axiom5}, and \eqref{axiom10}.  The other $5$ groups of conditions correspond to the remaining $5$ axioms in \eqref{eq:axioms}.
\end{proof}

\subsection{Filtration on $UA$}
\label{subsec:filtration}

Each homogeneous element in $UA$ has a length: Elements in $K$ has length $0$.  The elements $\gamma_*(a)$, where $\gamma = \alpha, \beta$, $* = l, r, m$, and $a \in A$, have length $1$.  Inductively, the homogeneous elements in $UA$ of length at most $k + 1$ ($k \geq 1$) are the $K$-linear combinations of the elements $\gamma_*(a) \cdot x$ and $x \cdot \gamma_*(a)$, where $x$ has length at most $k$.

For $k \geq 0$, consider the following submodule of $UA$:
   \[
   F_k UA = \lbrace x \in UA \colon x \text{ has length} \leq k \rbrace.
   \]
These submodules form an increasing and exhaustive filtration of $UA$, so that $F_0UA = K$, $F_kUA \subseteq F_{k+1}UA$, and $UA = \cup_{k\geq 0} F_kUA$.  Moreover, it is multiplicative, in the sense that $(F_kUA)\cdot(F_lUA) \subseteq F_{k+l}UA$.  Therefore, it makes sense to consider the associated graded algebra $Gr_* UA = \oplus_{k\geq 0} \, Gr_kUA$, where $Gr_k UA = F_kUA/F_{k-1}UA$.  It is clear that $Gr_0UA = K$.  The following result identifies the other associated quotients.

\begin{thm}
\label{thm:filtration}
In the associated graded algebra $Gr_* UA$, one always has $Gr_n UA = 0$ for all $n \geq 3$.  Moreover, there are isomorphisms
   \[
   \begin{split}
   Gr_1 UA ~\cong~ & \alpha_l A \oplus \alpha_r A \oplus \alpha_m A \oplus \beta_l A \oplus \beta_r A \oplus \beta_m A, \\
   Gr_2 UA ~\cong~ & (\alpha_l\beta_l A^{\otimes 2}) \oplus (\alpha_r\beta_r A^{\otimes 2}) \oplus (\alpha_r\beta_l A^{\otimes 2})  \\
   & \oplus (\alpha_m\beta_l A^{\otimes 2}) \oplus (\alpha_m\beta_r A^{\otimes 2}) \oplus (\alpha_r\beta_m A^{\otimes 2}) \oplus (\alpha_m\beta_m A^{\otimes 2}).
   \end{split}
   \]
In particular, if $A$ is of finite dimension $d$ over $K$, then $Gr_* UA$ has dimension $1 + 6d + 7d^2$ over $K$.
\end{thm}

\begin{proof}
The first isomorphism is clear, since $Gr_0 UA = F_0 UA = K$ and the elements of length $1$ are linearly generated by the $\alpha_*(a)$ and $\beta_*(a)$.  For the second isomorphism, observe that any other generator of length $\leq 2$ is identified with an element of length $1$ by one of the $33$ conditions, leaving only the $7$ displayed generators.  For example, $\beta_l(b) \cdot \beta_l(a) \in \beta_l\beta_l A^{\otimes 2}$ is identified with $\beta_l(a \lprod b) \in \beta_l A$ by condition (1) in the definition of the universal enveloping algebra $UA$.

Finally, to show that $Gr_n UA = 0$ for $n \geq 3$, it suffices to show that $F_2 UA = F_3 UA$.  It is straightforward to check that multiplying any one of the $7$ generators of length $2$ with a generator of length $1$ always yields an element of length $\leq 2$.  For example,
   \[
   \alpha_l(a) \cdot (\alpha_l(b) \cdot \beta_l(c))
   \,=\, (\alpha_l(a) \cdot \alpha_l(b)) \cdot \beta_l(c))
   \,=\, \alpha_l(a \lprod b) \cdot \beta_l(c),
   \]
which lies in $\alpha_l \beta_l A^{\otimes 2}$, by condition (7) in the definition of $UA$.
\end{proof}

An immediate consequence of this result is:

\begin{cor}[Poincar\'{e}-Birkhoff-Witt]
\label{cor:PBW}
Let $V$ denote the underlying vector space of a triassociative algebra $A$.  Equip $V$ with the abelian triassociative algebra structure ($a \ast b = 0$ always).  Then there exists an isomorphism $Gr_* UA \,\cong\, UV$ of unital associative algebras.
\end{cor}

\subsection{Corepresentation}
\label{subsec:corep}

By a \emph{corepresentation of} $A$, or an $A$-corepresentation, we mean a right $UA$-module.  This definition can be made more explicit as follows.  Let $N$ be an $A$-corepresentation.  Set
   \[
   \begin{cases}
   a < x := x \cdot \alpha_l(a) & \\
   a > x := x \cdot \alpha_r(a) & \\
   a \wedge x := x \cdot \alpha_m(a), 
   \end{cases}
   \begin{cases}
   x < a := x \cdot \beta_l(a) & \\
   x > a := x \cdot \beta_r(a) & \\
   x \wedge a := x \cdot \beta_m(a)
   \end{cases}
   \]
for $x \in N$ and $a \in A$.  This gives rise to three left actions $<, >, \wedge \colon A \otimes N \to N$ and three right actions $<, >, \wedge \colon N \otimes A \to N$.  With these notations, the condition that $N$ is an $A$-corepresentation is equivalent to the following $33$ axioms (for $x \in N$ and $a, b \in A$), corresponding to the axioms for $UA$ in \S \ref{subsec:uni}:
   \[
   \begin{split}
   &
   \begin{cases}
   (x < b) < a = x < (a \lprod b) = x < (a \,\mprod\, b) = x < (a \rprod b) & \\
   a < (x < b) = (a < x) < b = (a < x) \wedge b = (a < x) > b & \\
   (a \lprod b) < x = b < (a < x) = b \wedge (a < x) = b > (a < x) &
   \end{cases} \\
   &
   \begin{cases}
   x > (a \rprod b) = (x > b) > a = (x > b) < a = (x > b) \wedge a & \\
   (a > x) > b = a > (x > b) = a < (x > b) = a \wedge (x > b) & \\
   b > (a > x) = (a \rprod b) > x = (a \lprod b) > x = (a \,\mprod\, b) > x 
   \end{cases}   \\
   &
   \begin{cases}
   (x < b) > a = x > (a \lprod b) & \\
   a > (x < b) = (a > x) < b & \\
   (a \rprod b) < x = b < (a > x) 
   \end{cases}
   \begin{cases}
   (x < b) \wedge a = x \wedge (a \lprod b) & \\
   a \wedge (x < b) = (a \wedge x) < b & \\
   (a \,\mprod\, b) < x = b < (a \wedge x)
   \end{cases} \\
   & 
   \begin{cases}
   (x \wedge b) < a = x \wedge (a \rprod b) & \\
   a < (x \wedge b) = (a \wedge x) > b & \\
   (a \lprod b) \wedge x = b > (a \wedge x)
   \end{cases}
   \begin{cases}
   (x \wedge b) > a = x > (a \,\mprod\, b) & \\
   a > (x \wedge b) = (a > x) \wedge b & \\
   (a \rprod b) \wedge x = b \wedge (a > x)
   \end{cases} 
   \end{split}
   \]
   \[
   \begin{cases}
   (x \wedge b) \wedge a = x \wedge (a \,\mprod\, b) & \\
   (a \,\mprod\, b) \wedge x = b \wedge (a \wedge x) & \\
   a \wedge (x \wedge b) = (a \wedge x) \wedge b.
   \end{cases}   
   \]

Here are some examples of corepresentations:
   \begin{enumerate}
   \item $UA$ is an $A$-corepresentation via the right action of $UA$ on itself.
   \item $K$ is an $A$-corepresentation, where the $\alpha_*(a)$ and $\beta_*(a)$ in $UA$ act trivially.
   \item Let $M$ be an $A$-representation (= left $UA$-module).  It is straightforward to check that the definitions,
      \[
      \begin{cases}
      x \cdot \beta_l(a) := a \rprod x & \\
      x \cdot \alpha_r(a) := x \lprod a & \\
      x \cdot \beta_m(a) := a \,\mprod\, x & \\
      x \cdot \alpha_m(a) := x \,\mprod\, a,
      \end{cases}
      \begin{cases}
      x \cdot \beta_r(a) := 0 & \\
      x \cdot \alpha_l(a) := 0,
      \end{cases}
      \]
   for $x \in M$ and $a \in A$, give the underlying vector space of $M$ a right $UA$-module structure, denoted by $M^{op}$.  We call it the \emph{opposite corepresentation of} $M$.  In particular, considering $A$ as an $A$-representation, we have the opposite corepresentation $A^{op}$ of $A$.
   \end{enumerate}

Now suppose that $N$ is an $A$-corepresentation.  Define the \emph{module of $n$-chains of $A$ with coefficients in $N$} to be
   \[
   \CT_n(A,N) := K\lbrack T_n \rbrack \otimes N \otimes A^{\otimes n}.
   \]
Define a map
   \[
   d = \sum_{i=0}^n (-1)^i d_i \colon \CT_n(A,N) \to \CT_{n-1}(A,N),
   \]
where
   \[
   d_i(\psi \otimes x \otimes \mathbf{a}) = (d_i \psi) \otimes d_i^\psi(x \otimes \mathbf{a})
   \]
for $\psi \in T_n$, $x \in N$ and $\mathbf{a} = (a_1, \ldots , a_n) \in A^{\otimes n}$.  Here $d_i \psi$ is as in cohomology.  The maps $d^\psi_i$ are defined as follows.  Write $\psi$ uniquely as the grafting $\psi_0 \vee \cdots \vee \psi_k$ as before, where the valence of the lowest internal vertex in $\psi$ is $k + 1$.  Then set
   \[
   d^\psi_0(x \otimes \mathbf{a}) :=
   \begin{cases}
   (x \cdot \alpha_l(a_1)) \otimes (a_2, \ldots , a_n) & \text{if } \vert \psi_0 \vert = 0 \text{ and } k = 1 \\
   (x \cdot \alpha_m(a_1)) \otimes (a_2, \ldots , a_n) & \text{if } \vert \psi_0 \vert = 0 \text{ and } k > 1 \\
   (x \cdot \alpha_r(a_1)) \otimes (a_2, \ldots , a_n) & \text{if } \vert \psi_0 \vert > 0.
   \end{cases}
   \]
For $1 \leq i \leq n - 1$, set
   \[
   d^\psi_i(x \otimes \mathbf{a}) := x \otimes (a_1, \ldots , a_i \circ^\psi_i a_{i+1}, \ldots , a_n),
   \]
where $\circ^\psi_i$ is as in cohomology.  In particular, our $d_i$ coincides with the $d_i$ in Loday-Ronco \cite{lr} when $N = K$.  Finally, set
   \[
   d^\psi_n(x \otimes \mathbf{a}) :=
   \begin{cases}
   (x \cdot \beta_l(a_n)) \otimes (a_1, \ldots , a_{n-1}) & \text{if } \vert \psi_k \vert > 0 \\
   (x \cdot \beta_m(a_n)) \otimes (a_1, \ldots , a_{n-1}) & \text{if } k > 1 \text{ and } \vert \psi_k \vert = 0 \\
   (x \cdot \beta_r(a_n)) \otimes  (a_1, \ldots , a_{n-1}) & \text{if } k = 1 \text{ and } \vert \psi_1\vert = 0.
   \end{cases}
   \]
In order to show that $(\CT_*(A,N), d)$ is a chain complex, we make use of the cotangent complex.

\subsection{Cotangent complex}
\label{subsec:cotangent}

Consider $UA$ as an $A$-corepresentation.  Then the vector space $\CT_n(A,UA)$ becomes a left $UA$-module, where $UA$ acts only on the $UA$ factor.  Let $M$ be any $A$-representation ($=$ left $UA$-module).  There is an isomorphism
   \[
   \begin{split}
   \Hom_K(K\lbrack T_n \rbrack \otimes A^{\otimes n}, M) = \Ctrias^n(A,M) & ~\xrightarrow{\cong}~ \Hom_{UA}(\CT_n(A,UA), M) \\
   f & ~\mapsto~ (\overline{f} \colon \psi \otimes x \otimes \mathbf{a} \mapsto x \cdot f(\psi \otimes \mathbf{a})).
   \end{split}
   \]
It is clear that, under this isomorphism, the map $\delta^n_i \colon \Ctrias^n(A,M) \to \Ctrias^{n+1}(A,M)$ corresponds to $\Hom_{UA}(d_i,M)$.  Since $\Ctrias^*(A,M)$ is a cochain complex, we obtain:

\begin{prop}
\label{prop:cotangent}
$(\CT_*(A,UA),d)$ is a chain complex.
\end{prop}

Now observe that for an $A$-corepresentation $N$, there is an isomorphism
   \[
   \CT_n(A, UA) \, \xrightarrow[\cong]{N \otimes_{UA} (-)}\, \CT_n(A,N),
   \]
under which the $d_i$ in $\CT_*(A,UA)$ corresponds to the $d_i$ in $\CT_*(A,N)$.  Therefore, we conclude:

\begin{cor}
\label{cor:chain}
For any corepresentation $N$ of $A$, $(\CT_*(A,N),d)$ is a chain complex.
\end{cor}

\subsection{Triassociative homology}
\label{subsec:tri homology}

For an $A$-corepresentation $N$, define the \emph{triassociative algebra homology of $A$ with coefficients in $N$} as
   \[
   \HT_n(A,N) := H_n(\CT_*(A,N),d).
   \]

When $K$ is considered as an $A$-corepresentation, the module $\CT_n(A,K) \cong K \lbrack T_n \rbrack \otimes A^{\otimes n}$ is denoted by $\CT_n(A)$ in Loday-Ronco \cite{lr}.  Observe that the maps $d_0$ and $d_n \colon \CT_n(A,K) \to \CT_{n-1}(A,K)$ are both trivial.  In fact, for an element $\lambda \in K$,
   \[
   d_0^\psi(\lambda \otimes \mathbf{a}) = (\lambda \cdot \alpha_*(a_1)) \otimes (a_2, \ldots , a_n) = 0,
   \]
and similarly for $d_n^\psi$.  Therefore, in this case, the differential $d$ reduces to $d = \sum_{i=1}^{n-1} d_i$, which coincides with the differential $d \colon \CT_n(A) \to \CT_{n-1}(A)$ constructed in \cite{lr}.  In other words, we have
   \[
   \HT_n(A,K) = \HT_n(A),
   \]
where $\HT_n(A)$ is the triassociative algebra homology (with trivial coefficients) defined in \cite{lr}.

\subsection{$\HT_0$}
\label{subsec:H_0}
For example, identifying $\CT_1(A,N)$ with $N \otimes A$ and $\CT_0(A,N)$ with $N$, we have that $d = d_0 - d_1$, where
   \[
   \begin{cases}
   d_0(x \otimes a) & =\, x \cdot \alpha_l(a) \\
   d_1(x \otimes a) & =\, x \cdot \beta_r(a).
   \end{cases}
   \]
So the image of the differential $d \colon \CT_1(A,N) \to \CT_0(A,N)$ is the submodule $N \cdot (\alpha_l - \beta_r)$ of $N$, and, therefore, we have
   \[
   \HT_0(A,N) \,\cong\, \frac{N}{N \cdot (\alpha_l - \beta_r)}.
   \]
In particular:
   \begin{enumerate}
   \item If $N = K$, then $\HT_0(A,K) \cong K$.
   \item If $N = M^{op}$, the opposite corepresentation of an $A$-representation $M$, then $\HT_0(A,M^{op}) \cong M^{op}$.
   \end{enumerate}

\subsection{$\HT_1$}
\label{subsec:H_1}

Similarly, $\HT_1(A,N) = \ker d^1 /\im d^2$, where $\ker d^1$ is the submodule of $\CT_1(A,N) = N \otimes A$,
   \[
   \ker d^1 = \left\lbrace \sum_i x_i \otimes a_i \in N \otimes A ~ \colon ~ \sum_i a_i < x_i = \sum_i x_i > a_i\right\rbrace
   \]
The image of $d^2 \colon \CT_2(A,N) \to \CT_1(A,N) = N \otimes A$ is the span of the following three types of elements ($x \in N$, $a, b \in A$):
   \[
   \begin{split}
   & (a < x) \otimes b - x \otimes (a \lprod b) + (x < b) \otimes a, \\
   & (a > x) \otimes b - x \otimes (a \rprod b) + (x > b) \otimes a, \\
   & (a \wedge x) \otimes b - x \otimes (a \,\mprod\, b) + (x \wedge b) \otimes a.
   \end{split}
   \]
In particular, when $N = K$, we have 
   \[
   \HT_1(A,K) \cong \frac{A}{\text{span} \lbrace a \ast b \colon a, b \in A,~ \ast \,=\, \lprod, \rprod, \, \mprod \rbrace},
   \]
the abelianization of $A$.


\section{Deformations of triassociative algebras}
\label{sec:def}

In this final section, we describe algebraic deformations of a triassociative algebra $A = (A, \lprod, \rprod, \mprod)$, using the cohomology theory constructed in \S \ref{sec:def complex}.  The theory of triassociative algebra deformations is what one would expect from the existing literature on deformations.  The arguments in this section are rather straightforward and are very similar to the cases of associative algebras \cite{ger} and dialgebras \cite{mm}, whose arguments are provided in details and can be easily adapted to the present case.  Also, Balavoine \cite{bal} has another way of doing algebraic deformations that can be applied to triassociative algebras.   Therefore, we will safely omit most of the proofs in this section.

Thinking of $A$ also as an $A$-representation via the identity self-map, we consider the cochain complex $\Ctrias^*(A,A)$.  A $2$-cochain $\theta \in \Ctrias^2(A,A)$ will be identified with the triple $(\lambda, \rho, \mu)$, where each component is a binary operation on $A$, via
   \[
   \theta(\psi ; x, y) \,=\,
   \begin{cases}
   \lambda(x,y)& \text{if } \psi = \begin{picture}(16,10)     
                 \drawline(8,0)(8,2)(0,10)
                 \drawline(8,2)(16,10)
                 \drawline(12,6)(8,10)
                 \end{picture}, \\
   \rho(x,y)   & \text{if } \psi = \begin{picture}(16,10)     
                 \drawline(8,0)(8,2)(0,10)
                 \drawline(4,6)(8,10)
                 \drawline(8,2)(16,10)
                 \end{picture}, \\
   \mu(x,y)    & \text{if } \psi = \begin{picture}(16,10)     
                 \drawline(8,0)(8,10)
                 \drawline(8,2)(0,10)
                 \drawline(8,2)(16,10)
                 \end{picture},
   \end{cases}
   \]
for $x, y \in A$.

\subsection{Deformation and equivalence}
\label{subsec:def and eq}

By a \emph{deformation of $A$}, we mean a power series $\Theta_t = \sum_{i=0}^\infty \theta_i t^i$, in which each $\theta_i = (\lambda_i, \rho_i, \mu_i)$ is a $2$-cochain in $\Ctrias^2(A,A)$ and $\theta_0 = (\lprod, \rprod, \mprod)$, satisfying the following $11$ conditions (same as \eqref{eq:axioms}) for all $x, y, z \in A$, where $L_t = \sum_{i=0}^\infty \lambda_i t^i$, $R_t = \sum_{i=0}^\infty \rho_i t^i$, and $M_t = \sum_{i=0}^\infty \mu_i t^i$: (1) $L_t(L_t(x,y),z) = L_t(x, L_t(y,z))$, (2) $L_t(L_t(x,y),z) = L_t(x, R_t(y,z))$, etc.

Extending linearly, this gives a triassociative algebra structure on the power series $A \lbrack \lbrack t \rbrack \rbrack$.  We will also denote a deformation $\Theta_t$ by the triple $(L_t, R_t, M_t)$ of power series.

A \emph{formal isomorphism of $A$} is a power series $\Phi_t = \sum_{i=0}^\infty \phi_i t^i$, in which each $\phi_i \in \Ctrias^1(A,A) \cong \Hom_K(A,A)$ and $\phi_0 = \Id_A$.  Let $\Thetabar_t = (\Lbar_t, \Rbar_t, \Mbar_t)$ be another deformation of $A$.  Then $\Theta_t$ and $\Thetabar_t$ are said to be \emph{equivalent} if there exists a formal isomorphism $\Phi_t$ such that $\overline{O}_t(x,y) = \Phi_t O_t (\Phi^{-1}_t(x), \Phi^{-1}_t(y))$ for $O = L, R, M$ and for all $x, y, z \in A$.  In this case, we write $\Thetabar_t = \Phi_t \Theta_t \Phi_t^{-1}$.  Conversely, given a deformation $\Theta_t$ and a formal isomorphism $\Phi_t$, one can define a new deformation as $\Thetabar_t := \Phi_t \Theta_t \Phi_t^{-1}$.

\subsection{Infinitesimal}
\label{subsec:inf}

Given a deformation $\Theta_t = \sum_{i=0}^\infty \theta_i t^i$ of $A$, the linear coefficient $\theta_1$ is called the \emph{infinitesimal} of $\Theta_t$.

\begin{thm}
\label{thm:inf}
Let $\Theta_t = \sum_{i=0}^\infty \theta_i t^i$ be a deformation of $A$.  Then the infinitesimal $\theta_1$ is a $2$-cocycle in $\Ctrias^2(A,A)$, whose cohomology class is well-defined by the equivalence class of $\Theta_t$.  Moreover, if $\theta_i = 0$ for $1 \leq i \leq l$, then $\theta_{l+1}$ is a $2$-cocycle in $\Ctrias^2(A,A)$.
\end{thm}

This Theorem allows one to think of an infinitesimal as a cohomology class, instead of just a cochain.

\subsection{Rigidity}
\label{subsec:rigidity}

By the \emph{trivial deformation of $A$}, we mean the deformation $\Theta_t = (\lprod, \rprod, \mprod)t^0$.  We say that $A$ is  \emph{rigid} if every deformation of $A$ is equivalent to the trivial deformation.  The following result will lead to a cohomological criterion for rigidity.

\begin{prop}
\label{prop:rigid}
Let $\Theta_t = \theta_0 + \theta_l t^l + \theta_{l+1}t^{l+1} + \cdots$, with $\theta_i = (\lambda_i, \rho_i, \mu_i)$, be a deformation of $A$ for some $l \geq 1$ in which $\theta_l$ is a $2$-coboundary in $\Ctrias^2(A,A)$.  Then there exists a formal isomorphism of the form $\Phi_t = \Id_A + \phi t^l$ such that the deformation defined by $\Thetabar_t = \sum_{i=0}^\infty \thetabar_i t^i := \Phi_t \Theta_t \Phi_t^{-1}$ satisfies $\thetabar_i = 0$ for $1 \leq i \leq l$.
\end{prop}

Combining Theorem \ref{thm:inf} and Proposition \ref{prop:rigid}, we obtain the following cohomological criterion for rigidity.

\begin{cor}
\label{cor:rigid}
If $\Htrias^2(A,A)$ is trivial, then $A$ is rigid.
\end{cor}

Using Theorem \ref{thm:trivial Hi} and Corollary \ref{cor:rigid}, we obtain a large class of rigid objects.

\begin{cor}
\label{cor:free}
Free triassociative algebras are rigid.
\end{cor}

Next we want to obtain a cohomological criterion for the existence of a deformation with a prescribed $2$-cocycle as its infinitesimal.

\subsection{Deformations of finite order}
\label{subsec:finite}

Suppose that $N \in \lbrace 1, 2, \ldots  \rbrace$.  By a \emph{deformation of order $N$ of $A$}, we mean a polynomial $\Theta_t = \sum_{i=0}^N \theta_i t^i$, with each $\theta_i = (\lambda_i, \rho_i, \mu_i) \in \Ctrias^2(A,A)$ and $\theta_0 = (\lprod, \rprod, \mprod)$, satisfying the $11$ conditions for a deformation modulo $t^{N+1}$.

If $\theta_1 \in \Ctrias^2(A,A)$ is a $2$-cocycle, then $\theta_0 + \theta_1 t$ is a deformation of order $1$.  Thus, given a $2$-cocycle $\theta_1$, in order to determine the existence of a deformation with $\theta_1$ as its infinitesimal, it suffices to determine the obstruction to extending a deformation of order $N$ to one of order $N + 1$ for $N \geq 1$.

\subsection{Obstructions}
\label{subsec:obstructions}

Fix a deformation $\Theta_t = \sum_{i=0}^N \theta_i t^i$, with $\theta_i = (\lambda_i, \rho_i, \mu_i)$, of order $N < \infty$ as above.  Define a $3$-cochain $\Ob_\Theta \in \Ctrias^3(A,A)$ as follows.  Let $x, y, z$ be elements of $A$ and $\psi$ be a tree in $T_3$.  Then $\Ob_\Theta(\psi ; x, y, z)$ is the element (where $\Sigma = \Sigma_{j=1}^N$):
   \begin{subequations}
   \label{eq:Ob}
   \allowdisplaybreaks
   \begin{align}
   \sum \biglpren \lambda_j(\lambda_{N+1-j}(x,y),z) - \lambda_j(x, \lambda_{N+1-j}(y,z)) \bigrpren & \text{ if } \psi = \begin{picture}(18,11)    
                 \drawline(9,0)(9,2)(0,11)
                 \drawline(9,2)(18,11)
                 \drawline(12,5)(6,11)
                 \drawline(15,8)(12,11)
                 \end{picture}, \label{eq:Ob1} \\
   \sum \biglpren \lambda_j(\lambda_{N+1-j}(x,y),z) - \lambda_j(x, \rho_{N+1-j}(y,z)) \bigrpren & \text{ if } \psi = \begin{picture}(18,11)    
                 \drawline(9,0)(9,2)(0,11)
                 \drawline(9,2)(18,11)
                 \drawline(12,5)(6,11)
                 \drawline(9,8)(12,11)
                 \end{picture}, \label{eq:Ob2} \\
   \sum \biglpren \lambda_j(\rho_{N+1-j}(x,y),z) - \rho_j(x, \lambda_{N+1-j}(y,z)) \bigrpren & \text{ if } \psi = \begin{picture}(18,11)    
                 \drawline(9,0)(9,2)(0,11)
                 \drawline(3,8)(6,11)
                 \drawline(9,2)(18,11)
                 \drawline(15,8)(12,11)
                 \end{picture}, \label{eq:Ob3} \\
   \sum \biglpren \rho_j(\lambda_{N+1-j}(x,y),z) - \rho_j(x, \rho_{N+1-j}(y,z)) \bigrpren & \text{ if } \psi = \begin{picture}(18,11)    
                 \drawline(9,0)(9,2)(0,11)
                 \drawline(6,5)(12,11)
                 \drawline(9,8)(6,11)
                 \drawline(9,2)(18,11)
                 \end{picture}, \label{eq:Ob4} \\
   \sum \biglpren \rho_j(\rho_{N+1-j}(x,y),z) - \rho_j(x, \rho_{N+1-j}(y,z)) \bigrpren & \text{ if } \psi = \begin{picture}(18,11)    
                 \drawline(9,0)(9,2)(0,11)
                 \drawline(3,8)(6,11)
                 \drawline(6,5)(12,11)
                 \drawline(9,2)(18,11)
                 \end{picture}, \label{eq:Ob5} \\
   \sum \biglpren \lambda_j(\lambda_{N+1-j}(x,y),z) - \lambda_j(x, \mu_{N+1-j}(y,z)) \bigrpren & \text{ if } \psi = \begin{picture}(18,11)    
                 \drawline(9,0)(9,2)(0,11)
                 \drawline(9,2)(18,11)
                 \drawline(12,5)(6,11)
                 \drawline(12,5)(12,11)
                 \end{picture}, \label{eq:Ob6} \\
   \sum \biglpren \lambda_j(\mu_{N+1-j}(x,y),z) - \mu_j(x, \lambda_{N+1-j}(y,z)) \bigrpren & \text{ if } \psi = \begin{picture}(18,11)    
                 \drawline(9,0)(9,2)(0,11)
                 \drawline(9,2)(18,11)
                 \drawline(9,2)(9,11)
                 \drawline(15,8)(12,11)
                 \end{picture}, \label{eq:Ob7} \\
   \sum \biglpren \mu_j(\lambda_{N+1-j}(x,y),z) - \mu_j(x, \rho_{N+1-j}(y,z)) \bigrpren & \text{ if } \psi = \begin{picture}(18,11)    
                 \drawline(9,0)(9,2)(0,11)
                 \drawline(9,2)(18,11)
                 \drawline(9,2)(9,8)
                 \drawline(9,8)(6,11)
                 \drawline(9,8)(12,11)
                 \end{picture}, \label{eq:Ob8} \\
   \sum \biglpren \mu_j(\rho_{N+1-j}(x,y),z) - \rho_j(x, \mu_{N+1-j}(y,z)) \bigrpren & \text{ if } \psi = \begin{picture}(18,11)    
                 \drawline(9,0)(9,2)(0,11)
                 \drawline(9,2)(18,11)
                 \drawline(9,2)(9,11)
                 \drawline(3,8)(6,11)
                 \end{picture}, \label{eq:Ob9} \\
   \sum \biglpren \rho_j(\mu_{N+1-j}(x,y),z) - \rho_j(x, \rho_{N+1-j}(y,z)) \bigrpren & \text{ if } \psi = \begin{picture}(18,11)    
                 \drawline(9,0)(9,2)(0,11)
                 \drawline(9,2)(18,11)
                 \drawline(6,5)(6,11)
                 \drawline(6,5)(12,11)
                 \end{picture}, \label{eq:Ob10} \\
   \sum \biglpren \mu_j(\mu_{N+1-j}(x,y),z) - \mu_j(x, \mu_{N+1-j}(y,z)) \bigrpren & \text{ if } \psi = \begin{picture}(18,11)    
                 \drawline(9,0)(9,2)(0,11)
                 \drawline(9,2)(18,11)
                 \drawline(9,2)(6,11)
                 \drawline(9,2)(12,11)
                 \end{picture}.
   \end{align}
   \end{subequations}
Notice that the definition of the $3$-cochain $\Ob_\Theta$ arises from the coefficients of $t^{N+1}$ in the $11$ conditions for a deformation by subtracting the right-hand side from the left-hand side and deleting the $4$ terms corresponding to $j = 0$ and $j = N + 1$.

\begin{lemma}
\label{lem:ob}
The element $\Ob_\Theta \in \Ctrias^3(A,A)$ is a $3$-cocycle.
\end{lemma}

The proof of the Lemma is a long but elementary computation.  The following result is independent of Lemma \ref{lem:ob}.

\begin{thm}
\label{thm:ob}
Let $\Theta_t = \sum_{i=0}^N \theta_i t^i$, with $\theta_i = (\lambda_i, \rho_i, \mu_i)$, be a deformation of order $N < \infty$ of $A$ and $\theta_{N+1}$ be an element of $\Ctrias^2(A,A)$.  Define the polynomial $\Thetatilde_t := \Theta_t + \theta_{N+1}t^{N+1}$.  Then $\Thetatilde_t$ is a deformation of order $N + 1$ if and only if $\Ob_\Theta = \delta^2 \theta_{N+1}$.
\end{thm}

Combining Theorem \ref{thm:ob} with Lemma \ref{lem:ob}, one concludes that the obstruction to extending a deformation $\Theta_t$ of order $N$ to one of order $N + 1$ is the cohomology class of $\Ob_\Theta$, which lies in $\Htrias^3(A,A)$.  Therefore, all these obstructions vanish if the cohomology module $\Htrias^3(A,A)$ is trivial.

\begin{cor}
\label{cor:ob}
If $\Htrias^3(A,A)$ is trivial, then every $2$-cocycle in $\Ctrias^2(A,A)$ is the infinitesimal of some deformation of $A$.
\end{cor}



\begin{thebibliography}{99}

\bibitem{bal}D.\ Balavoine, Deformations of algebras over a quadratic operad, Contemp.\ Math.\ \textbf{202} (1997), 207-234.

\bibitem{bok}M.\ B\"{o}kstedt, Topological Hochschild homology, University of Bielefeld preprint, 1986.

\bibitem{ekmm}A.\ D.\ Elmendorf, I.\ Kriz, M.\ A.\ Mandell and J.\ P.\ May, Rings, modules, and algebras in stable homotopy theory, with an appendix by M.\ Cole.  Math.\ Surveys and Monographs, vol. \textbf{47}, Amer.\ Math.\ Soc., Providence, RI, 1997.

\bibitem{frabetti}A.\ Frabetti, Dialgebra (co)homology with coefficients, in: Dialgebras and related operads, 67-103, Lecture Notes in Math. \textbf{1763}, Springer, Berlin, 2001.

\bibitem{ger}M.\ Gerstenhaber, On the deformation of rings and algebras, Ann.\ Math.\ \textbf{79} (1964), 59-103.

\bibitem{gk}V.\ Ginzburg and M.\ M.\ Kapranov, Koszul duality for operads, Duke Math.\ J.\ \textbf{76} (1994), 203-272.

\bibitem{loday}J.-L.\ Loday, Dialgebras, in:  Dialgebras and related operads, 7-66, Lecture Notes in Math. \textbf{1763}, Springer, Berlin, 2001.

\bibitem{lr}J.-L.\ Loday and M.\ Ronco, Trialgebras and families of polytopes, in: Homotopy theory: relations with algebraic geometry, group cohomology, and algebraic $K$-theory, 369-398, Contemp.\ Math. \textbf{346}, Amer.\ Math.\ Soc., Providence, RI, 2004.

\bibitem{mm}A.\ Majumdar and G.\ Mukherjee, Deformation theory of dialgebras, K-theory \textbf{27} (2002), 33-60.

\bibitem{ms}J.\ E.\ McClure and J.\ H.\ Smith, A solution of Deligne's Hochschild cohomology conjecture, 153-193, Contemp.\ Math.\ \textbf{293}, Amer.\ Math.\ Soc., Providence, RI, 2002.

\bibitem{sta}J.\ D.\ Stasheff, Homotopy associativity of $H$-spaces. I, II.  Trans.\ Amer.\ Math.\ Soc.\ \textbf{108} (1963), 275-292 and 293-312.

\bibitem{weibel}C.\ A.\ Weibel, An introduction to homological algebra, Cambridge studies in advanced mathematics \textbf{38}, Cambridge Univ.\ Press, Cambridge, UK, 1994.

\end{thebibliography}
\end{document}